\documentclass[conference]{IEEEtran}
\IEEEoverridecommandlockouts
\usepackage{amsmath,amssymb,amsfonts}
\usepackage{algorithmic}
\usepackage{cite}
\usepackage{graphicx}
\usepackage{framed}
\usepackage{mathtools}
\usepackage{subcaption}
\usepackage{textcomp}
\usepackage[normalem]{ulem} 
\usepackage{xcolor}
\usepackage{hyperref}

\newcommand{\expect}{\operatorname{\mathbb{E}}}
\newcommand{\prob}{\operatorname{\mathbb{P}}}

\begin{document}

\title{
    Baseline hydropower generation offer curves 
}

\author{
 \IEEEauthorblockN{Jon Pearce, Arash Khojaste, Golbon Zakeri}
    \IEEEauthorblockA{\textit{Mechanical and Industrial Engineering} \\
        \textit{The University of Massachusetts Amherst}\\
        Amherst, MA \\
        \{jonathanpear, akhojaste, gzakeri\}@umass.edu
    }
    \and
    \IEEEauthorblockN{Geoffrey Pritchard}
    \IEEEauthorblockA{\textit{Department of Statistics} \\
        \textit{The University of Auckland}\\
        Auckland, New Zealand \\
        https://orcid.org/0000-0002-5660-4503}
}

\maketitle

\begin{abstract}
    We outline a mathematical model for pricing hydropower generation.
    The model involves a Markov decision process that reflects the seasonal variation in historical time series of water inflows. 
    The procedure is computationally efficient and easy to interpret.
\end{abstract}


\section*{Nomenclature}

\subsection*{Sets and their elements}
\begin{tabular}{p{35pt}p{300pt}}
    \(r \in R\) & Inflow regimes (quantile sub-intervals)\\
    \(t \in T\) & Time points\\
    \(f \in F\) & Inflows\\
    \(f^\prime\) & Successor of inflow \(f\)
    \\
    \(\ell \in L\) & Water level in the reservoir\\
    \(s \in S\) & System states, \(S = L \times R \times T\)
    \\
    \(r^\prime \in R\) & Successor of regime \(r\) 
    \\
    \(s^\prime \in S\) & Successor of state \(s\)
    \\
    \(a \in A\) & Actions, reservoir outflow 
\end{tabular}   

\subsection*{Fourier regression models}
\begin{tabular}{p{35pt}p{300pt}}
    \(\vec{\phi}(t)\) & Fourier basis functions 
    \\
    \(\omega\)         & Frequency in trigonometric polynomials \\
    \(\vec{\beta}, \vec{\gamma}\) & Regression parameters
    \end{tabular} 
 
\subsection*{Inflow Markov chain}
\begin{tabular}{p{35pt}p{300pt}}
    \(p_{r^\prime|rt}\) & Probability of transition \(r \rightarrow r^\prime\) at time \(t\)
    \\
    \(p_{f|rt}\) & Probability of inflow \(f\) in regime \(r\) at time \(t\)
    \\
    \(\sim\) & Sampling from a distribution
\end{tabular} 

\subsection*{Operating Markov decision model}
\begin{tabular}{p{35pt}p{300pt}}
    \(p_{s^\prime|sa}\) & Probability of transition \(s \rightarrow s^\prime\) under \(a\)
    \\
    \(c_{sa}\) & Instantaneous cost of action \(a\) in state \(s\)
    \\
    \(y^\ast_{sa}\)& state-action distribution of water release policy
\end{tabular}    

\subsection*{Pricing Markov model}
\begin{tabular}{p{35pt}p{300pt}}
    \(u\) & Expected annual cost \\
    \(v^\ast_s\)& Value to the system of state \(s\)
\end{tabular}    

\section{Introduction}
Many countries and jurisdictions around the world, including most of the United States, meet their electricity needs through organized electricity markets. In these markets, generators submit supply offers that reflect their marginal cost of production. Calculating this cost is relatively straightforward for traditional power sources, such as fossil fuel-fired generation units. For example, natural gas turbines can base their offers on known fuel contract prices and turbine power curves, which describe power output as a function of fuel input. Additionally, the costs of operations and maintenance can be estimated with reasonable accuracy at the time of dispatch.

In contrast, estimating the marginal cost of electricity production from hydropower resources is substantially more complex. Offer stacks for hydro reservoirs must reflect the trade-off between immediate power production and the value of conserving water for future use (the so-called opportunity cost), especially given uncertainties in future inflows, fluctuating demand, and volatile electricity prices. Accurately estimating and incorporating this water value into the electricity market is critical for balancing the flexibility of thermal generation against the strategic use of limited hydro resources, ultimately enhancing both reliability and economic efficiency. A reasonable estimate of water values can also serve as a benchmark to assess whether a generator is acting competitively in the market or exerting market power. Furthermore, as renewable generation sources such as wind and solar achieve greater penetration, and as demand and transportation sectors become increasingly electrified, the overall system becomes more stochastic and volatile. This heightened uncertainty further complicates the task of valuing water accurately, making water value estimation increasingly important for efficient market participation and system planning.

One possible way to evaluate water over a time horizon is to use stochastic dynamic programming (SDP). This technique explicitly models the sequential and uncertain nature of water inflows and operational decisions over time. By breaking the problem into stages and optimizing expected future profits at each stage, SDP helps identify water release and storage strategies that balance immediate generation with future value. Pereira and Pinto \cite{Pereira-Pinto:1991} introduced Stochastic Dual Dynamic Programming (SDDP) as a scalable method to value hydropower systems under uncertainty. Their approach uses dual information to build piecewise linear approximations of future cost functions, enabling efficient planning of water use across large, complex hydro systems with uncertain inflows. In order to utilize SDDP-like approaches, we need to construct a scenario tree to model inflows over the time horizon, which can prove an arduous task. Furthermore, convergence of SDDP can stall; for instance, gaps of nearly 22\% are reported for an instance of the long-term planning problem of the Brazilian power system \cite{Shapiro-et-al:2013}.

In what follows, we will describe a Markov decision process approach to obtaining water values for a reservoir. Our methodology uses historical inflow data over decades to build a Markov model of inflow following \cite{Pritchard:2015:Inflow-Modeling}. To account for the periodic (seasonal) nature of inflow, we fit quantile Fourier regressions that capture the stochastic dynamics of inflows over time. Specifically, we develop a time-inhomogeneous Markov model that predicts the inflow volume (more precisely, its corresponding quantiles) for each week of the year, conditioned on the inflow observed in the preceding week. We then utilize this model to find a plan of operations of a reservoir that minimizes the total average long-run cost of the system. Coupled with such an optimized plan is a set of water value curves that provide the value of water in a reservoir for each week of the year. 

Our main contributions in this study include: 
\begin{itemize}
    \item Formulating a Markov decision process-based procedure for pricing hydropower generation, reflecting the seasonal variation in historical time series. The procedure is computationally efficient and provides a principled means of generating price scenarios over long-term investment planning horizons.
    \item Implementing this methodology on a realistic example based on New Zealand, utilizing 74 years of historical inflow data from that country (1948--2021). 
\end{itemize}

\section{Related work}

Catalao et al. (\cite{Catalao-et-al:2009:Nonlinear}, \cite{Catalao-et-al:2012:Optimal-Hydro}) propose a non-linear approach for scheduling releases in Portugal's reservoir cascade system that accounts for head effects and the size of the dam and they further develop the model to account for stochastic prices. Garcia-Gonzalez et al. \cite{GarciaGonzalez-et-al:2006:Generation-Bids} and Sauhats et al. (\cite{Sauhats-et-al:2015:Stochastic-Approach}) develop a stochastic optimization procedure to plan the operation of three hydropower plants. Their proposed stochastic optimization algorithm is based on time-average revenue maximization and utilizes neural networks. Garcia-Gonzalez et al. (\cite{GarciaGonzalez-et-al:2006:Generation-Bids}) develop a mixed-integer linear program that takes account of a detailed representation of the generating units. They maximize the expected profit of the generator responding to stochastic prices that are modeled using hidden Markov models. 

Kleiven et al. (\cite{Kleiven-Risanger-Fleten:2023:Co-movements}) develop a price model that considers the joint dynamics of forward prices and inflows to assess the frequently used assumption of independence of inflows from reservoir levels. They show that ignoring this dependence leads to underestimation of water values by the producers. Ilak et al. \cite{Ilak-et-al:2013:Profit-Maximization}, \cite{Ilak-et-al:2014:Price-Taker}) develop models for the self-scheduling problem of a hydro producer formulated and solved as a mixed integer linear programming problem. Their model is deterministic but accounts for head effects of the reservoir and considers participation in the ancillary services market. Daglish et al. (\cite{Daglish-et-al:2021:Pricing-Effects}) analyze the Brazilian forward market. The Brazilian electricity market is a hydro-dominated electricity market. They assessed generator entry before and after the forward market was introduced. Zhu et al. (\cite{Zhu-et-al:2006:Portfolio-MDP}) use MDPs to manage a portfolio optimization for hydro power contracts within a two-stage model. The operational model (in the second stage) is the conventional SDDP model. The above papers either respond to electricity market prices or focus on contract markets. In contrast, our approach is to find the true value of water/reservoir storage to the electricity system. 

Perhaps the most relevant comparison to our work is the work of \cite{Garcia-et-al:2001:Strategic}. That paper develops a simplified oligopoly model where hydro generators engage in dynamic Bertrand competition using a Markov strategy based on the state of water reservoirs at the beginning of each period. Our results provide guidance on what true water values, submitted by a {\em competitive} hydropower generator, should look like. This can aid in assessing whether a hydropower participant in the power market is acting as a competitive agent as opposed to exerting market power. 

\section{Model formulation}

The formulation that we describe provides an operating policy and associated benchmark prices for hydropower generation.

\subsection{Simplifications}

For simplicity, we formulate the model for a single reservoir; however, the model is easily generalizable to multiple reservoirs, and the results obtained in this study will (appropriately) extend to that case. The electrical load [MW] and cost of thermal generation [\$/MWh] are assumed constant when deriving the operating policy. When evaluating the policy through simulation, both may vary.


\subsection{A Markov chain model of reservoir inflow}\label{sec:inflow-mc}
 
The \(\alpha\)-quantiles \(q_\alpha(t)\) of inflow are assumed periodic and approximated from historical weekly time series via Fourier basis regression: 
\begin{align*}
    \omega = {}& 2\pi/365.25 \quad \text{rad/day}
    \\
    \vec{\phi_q}(t) = {}& [1,~ \cos(\omega t),~ \sin(\omega t),~ \cos(2 \omega t),~ \sin(2 \omega t)]
    \\
    q_\alpha(t|\vec{\beta}_\alpha) = {}& \vec{\phi_q}(t)\cdot\vec{\beta}_\alpha
\end{align*}
Regression coefficients \(\beta_\alpha\) are determined for a sequence of pre-selected levels \(0\% < \alpha_1 < \cdots < \alpha_{\vert R \vert - 1} < 100\%\) via quantile regression \cite{Khojaste-Pritchard-Zakeri:2024:QFR}, which entails solution of a linear program.

The time-varying inter-quantile intervals \([q_{\alpha_{r-1}}(t), q_{\alpha_r}(t)]\) define inflow regimes indexed by \(r \in R = \{1, 2, ...\}\), where \(q_{\alpha_0}(t) = 0\) and \(q_{\alpha_{|R|}} = \infty\).

The conditional distribution \(p_{f|rt} = \prob(f | r, t)\)--with support in \([q_{\alpha_{r-1}}(t), q_{\alpha_r}(t)]\)--of weekly inflow $f$ [MW], can be approximated from the inflow time series via histogram approximation.

Bridging these conditional distributions, the flow regime weekly series is modeled as a Markov chain with time-dependent, periodic state-transition probabilities:
\begin{align*}
    \vec{\phi_p}(t) = {}& [1,~ \cos(\omega t),~\sin(\omega t)]
    \\
    p_{r^\prime|rt}(\vec{\gamma}) = {}& \prob(r^\prime | r, t, \vec{\gamma}) = \vec{\phi_p}(t)\cdot\vec{\gamma}
\end{align*}
The regression parameters $\vec{\gamma}$--coefficients of the Fourier terms that provide the best fit to the historical regime sequence--are determined via maximization of the associated log-likelihood function, subject to (convex quadratic) non-negativity constraints and (linear) normalization constraints:
\begin{align*}
    \arg\max_{\vec{\gamma}} ~ \sum_{(t,r,r^\prime) \in \text{data}} \log p_{r^\prime|rt} (\vec{\gamma}) {}&
    \\
    \sum_{r^\prime} p_{r^\prime|rt}(\vec{\gamma}) = {}& 1 \quad \forall t, r
    \\
    p_{r^\prime|rt}(\vec{\gamma}) \in {}& [0, 1] \quad \forall t, r, r^\prime
\end{align*}
This problem reduces \cite{Khojaste-Pritchard-Zakeri:2024:QFR} to a second-order cone-constrained convex program \cite{Lobo-et-al:1998:SOCP}.

\subsection{Operational Markov decision model}

The operational characteristics of the reservoir can be expressed in terms of sets of states \(S\) and actions \(A\):
\begin{align}
    S = \{\ell, r, t | \ell \in L, r \in R, t \in T\}
\end{align}
where \(L\) and \(T\) are discrete approximations of the reservoir level (expressed in equivalent MW) and time (in weeks), while $A$ is a discrete approximation of the continuous set of operator-controlled outflows.

The time-varying inflow distributions described in \S\ref{sec:inflow-mc} induce a distribution \(p_{s^\prime|sa} = \prob(s^\prime | s, a)\) on state transitions \(s \overset{a}{\rightarrow} s^\prime\), as follows:
\begin{align*}
    \ell, r, t = s
    \\
    f \sim {}& p_{f|r t}
    \\
    r^\prime \sim {}& p_{r^\prime|r t}
    \\
    \ell^\prime = {}& \ell + f - a
    \\
    t^\prime = {}& t \operatorname{\%} 52 + 1
    \\
    s^\prime = {}& (\ell^\prime, r^\prime, t^\prime)
    \\
    p_{s^\prime|sa} = {}& \prob(s^\prime | s, a) = p_{f|rt} ~ p_{r^\prime|rt}
\end{align*}
 Note that energy (MWh) is expressed as a rate (MW) with a fixed time step of one week.
The formula for \(\ell^\prime\) is illustrative only: A real implementation must ensure that \(\ell\) is non-negative and within the reservoir's storage capacity.


In addition to a state transition probability, each state-action pair \(s-a\) is associated with a financial cost. Although water discharge, itself, has zero marginal cost, \(a\) determines the net demand that must be met with thermal generation (nonzero marginal cost) or by demand curtailment (punitive cost).

\begin{align*}
    \delta = {}& \text{demand [MW]}
    \\
    \tau = {}& \text{thermal capacity [MW]}
    \\
    c_{sa} = {}&
    \min(\delta - a, \tau) \times \text{(fuel price)}  +
    \\
    {}& 
    \max(\delta - \tau - a, 0) \times \text{(curtailment price)}
\end{align*}

With discretizations \(S\) and \(A\), the costs \(c_{sa}\) and transition probabilities \(p_{s^\prime|sa}\) constitute a Markov decision model, to be solved to determine a minimum-cost policy mapping each \(s \in S\) to an \(a \in A\).

\subsection{Minimum expected cost operating policy}


Denote by \(y_{sa}\) the probability of being in state \(s\) and executing action \(a\): The feasible set of the following LP \eqref{eqn:state-action-lp} comprises (unconditional) distributions over $S\times A$ that are consistent with the conditional probabilities \(p_{s^\prime|s a}\), and the objective function evaluates the expected total annual operating cost:
\begin{align}
    \begin{split}
        \vec{y}^\ast = \arg\min_{y \geq 0} ~ \sum_{s, a} c_{sa} y_{sa} &
        \\
        \sum_a y_{s^\prime a} = {}&
                \sum_{a, s} p_{s^\prime|sa} y_{sa} ~
            \forall s^\prime \in S
            \\
            \sum_{s, a} y_{sa} =  {}& 1
            \label{eqn:state-action-lp}
        \end{split}
\end{align}
Although \(\vec{y}^\ast\) ostensibly encodes a probabilistic operating policy (i.e., ``take action $a$ in state $s$ with probability $y^{\ast}_{sa}$''), it can be shown that, in the absence of floating-point error, \(y^\ast\) has at most one supported action for each state.

Since \eqref{eqn:state-action-lp} relies entirely on historical data, we should not expect \(\vec{y}^\ast\) to be competitive with alternative procedures that utilize online/rolling forecasts of inflow and demand data. However, it provides a benchmark for use in generation expansion planning and--as we describe next--provides benchmark offer stacks.

\subsection{Pricing formulation}

LP \eqref{eqn:state-action-lp} can be paired with a dual LP:
\begin{align}
    \begin{split}
        u^\ast, \vec{v}^\ast = 
        \arg\max_{u, \vec{v}} ~ u & 
        \\
         u + v_s \leq {}& c_{sa} + \sum_{s^\prime} p_{s^\prime|sa} v_{s^\prime} \quad  \forall s, a 
        \label{eq:val_max}
    \end{split}
\end{align}
By linear programming duality \cite{Boyd-Vandenberghe:2004:Convex, Ferris-Mangasarian-Wright:2007:LP-MATLAB}, \(u^\ast = \sum_{s, a} c_{sa} y^\ast_{sa}\) -- the maximum expected value or, equivalently, minimum expected cost of operating the integrated electricity system. 

The right-hand side of the constraint can be written \(c_{sa} + \expect(v_{s^\prime} | s, a)\), combining the immediate cost \(c_{sa}\) with an expectation over downstream states. This left-hand side involves another sum with cost \(u\), while \(\max u\) serves to tighten the inequality. Altogether, \(v^\ast_s\) provides a rational lower bound on the expected value of state \(s\), consistent with \(\{c_{sa}: a \in A\}\) and transitions \(s \rightarrow s^\prime\) under policy \(y^\ast\) (see, e.g., \cite{bertsekas_2018v2, Puterman:2005:MDP-Book}).

\section{Case study}\label{sec:case-study}

\begin{table}[]
    \centering
    \caption{Problem parameters}
    \begin{tabular}{rc}
        Description & Value \\ \hline
        reservoir storage capacity & 840\,GWh \\
        hydropower generation capacity & 900\,MW \\
        thermal generation capacity & 900\,MW \\
        load & 1400\,MW \\
        thermal fuel price & \$50/MWh \\
        load curtailment price & \$1000/MWh
    \end{tabular}
    \label{tab:parame}
\end{table}

\subsection{Context}

In this section, we apply our method to hydropower reservoir management to illustrate how we might obtain water value curves. 

The inflow data used in this study are the historical weekly inflow observed in the Waitaki catchment in New Zealand from 1948 to 2021 (see \cite{Pritchard:2015:Inflow-Modeling}). 
Although the Waitaki River system actually has multiple catchments and reservoirs, we simplify it here as a single equivalent reservoir. While the actual system has an energy storage capacity of about 2500\,GWh, we reduce it to 840\,GWh (\(\text{5000\,MW}\times \text{$(7 \times 24)$ h/week}\)) in our example to increase the likelihood of shortages and overflows and thus make the problem more illustrative.

Our system must meet a constant demand of 1400\,MW at the lowest possible expected cost. Supply consists of {900\,MW} of fossil-fueled generation (at \$50/MWh) and up to {1400\,MW} of hydropower generation, which is cost-free but subject to inflow uncertainty and limited storage capacity. Any unmet demand incurs a penalty cost of \$1000/MWh. The hydro consists of 500 MW equivalent of so called ``run of river'', plus an additional 900 MW that can be produced by running a turbine. The reservoir itself has a capacity of holding water that is equivalent to 5000 MW (which we represent as 50 segments of 100 MW worth of water).  

\subsection{Markov model}

\begin{figure}[!t]
  \centering
  \includegraphics[width=0.48\textwidth]{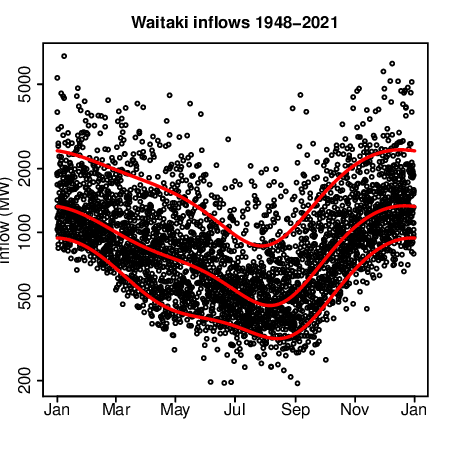}
  \caption{Weekly inflow data for the Waitaki system, with fitted 10th, 50th, and 90th
percentile models partitioning the data into four inflow states}
  \label{fig:WaitakiQFR}
\end{figure}

For this case study, we choose $$\alpha \in \{0\%, 10\%, 50\%, 90\%, 100\%\},$$
defining four inflow regimes \(r \in \text{1:4}\). Figure \ref{fig:WaitakiQFR} illustrates the quantile Fourier regression fits to the inflow data for the system of our case study. 

The reservoir inflows and outflows and non-hydro generation are binned into 100\,MW segments. 

We discretize the reservoir storage into 51 states (0:50, inclusive), each representing a block of 16.8\,GWh (equivalent to 100\,MW over a one week time step). Combined with the four inflow regimes and 52 weekly time points per operating period, the state set has size $$|S| = |L| \times |R|\times |T| = 51 \times 4 \times 52 = 10608.$$ 
    

At each time step, the available actions consist of releasing water for hydropower generation; since this must be delivered in 100\,MW blocks
there are \(|A| = 10\) different actions available. Thermal generation of 900\,MW is used to meet the remaining demand, with any shortfall resulting in curtailment.


Thus, the linear program \eqref{eqn:state-action-lp} (respectively, \eqref{eq:val_max}) has \(|S|\,|A| = 106080\) variables (constraints)  and \(|S| + 1 = 10608\) constraints (variables).

\section{Results}

\begin{figure*}[t]
    \centering
    \begin{subfigure}{0.6\textwidth}
        \includegraphics[width=\textwidth]{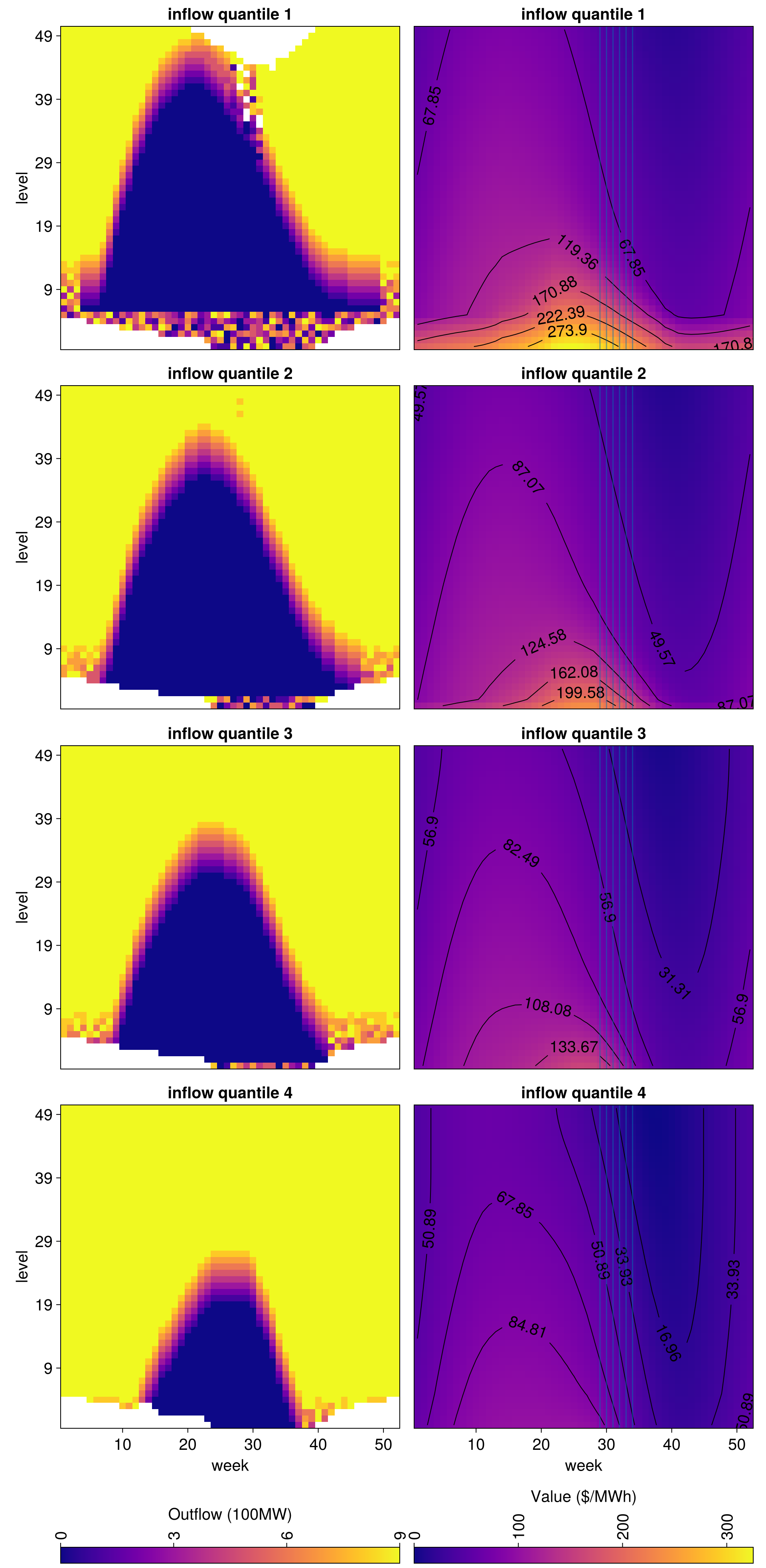}
        \caption{Optimal operating policy (left) and value function (right).}
        \label{fig:week-level}
    \end{subfigure}
    ~
    \begin{subfigure}{0.3\textwidth}
        \includegraphics[width=\textwidth]{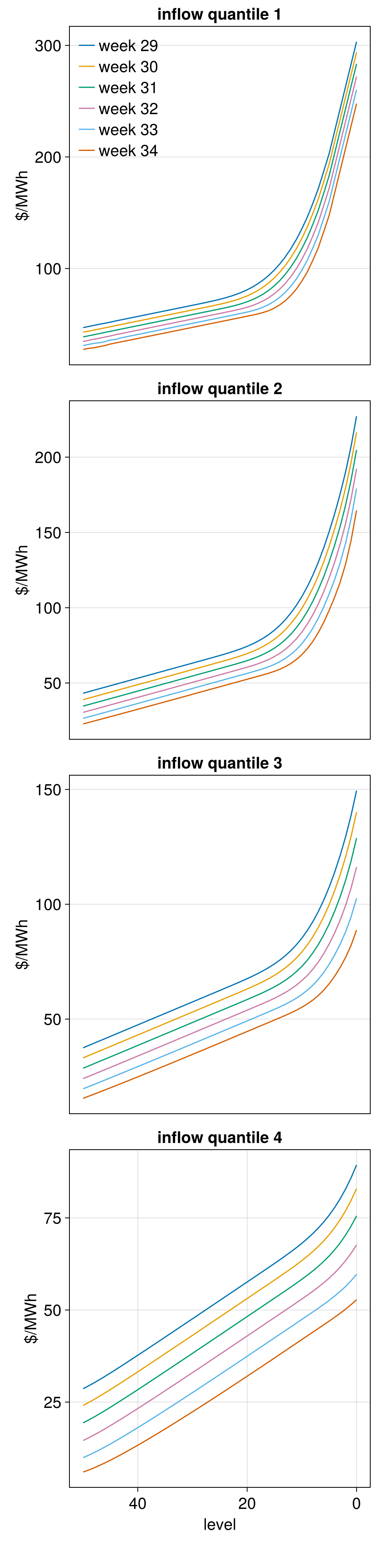}
        \caption{Cross sections of the value function.}
        \label{fig:level-price}
    \end{subfigure}
    \caption{Results of the case study in Section \ref{sec:case-study}. Vertical lines on the right of \ref{fig:week-level} show the sections of \ref{fig:level-price}.}\label{fig:results}
\end{figure*}

Figure \ref{fig:week-level} depicts the optimal solutions of the dual pair of LPs \eqref{eqn:state-action-lp} and \eqref{eq:val_max}. The left panel depicts the optimal water release policy as a function of the week and inflow regime, as extracted from \(\vec{y}_{sa}\). In New Zealand's drier winter weeks (25-38, say), the optimal policy withholds water in lower inflow regimes, even at high lake levels. For instance, discharge is zero in the lowest-inflow regime \(r = 1\) even when the reservoir is at 80\% capacity. Uncolored patches in this plot correspond to unsupported states, i.e., states that the model expects will never be occupied under the optimal policy. 

The middle panel shows the value to the system \(v^\ast_s\) of each reservoir state \(s\). As one might have anticipated, these values increase with decreasing storage level and in dryer inflow regimes.

The faint blue vertical lines in Figure \ref{fig:week-level} indicate cross-sections of the value surface, shown as curves in Figure \ref{fig:level-price}. 
Each cross-section defines an offer stack for water value. The value curves capture the expected value to the system of having a reservoir in various states (i.e., levels). This value (to the system) would reflect the fair payment that the system as a whole should be willing to pay to the hydro generator to supply electricity in accordance with the system's needs (and, equivalently, value the water highly and preserve its use for the future if the state of the reservoir and inflow so indicate). The water value curves are depicted in the right-most panel of the Figure \ref{fig:results}. Note that if we are looking to extract an estimated water value for a particular week, e.g., week 32 of the year, it is worthwhile having a distribution of water value curves around the ``trend" curve for week 32. This could be achieved by plotting curves for weeks similar to 32, namely a few weeks before and after week 32, as illustrated in our figure. 

\section{Conclusion and extensions}

We described a procedure for determining baseline offer curves for competitive hydropower generation, reflecting opportunity cost and inflow uncertainty.

Inflows were discretized using time-varying quantile intervals (regimes) via Fourier basis regression, and the regime dynamics were modeled as a Markov chain. The reservoir operational model comprises a state set of inflow regimes, storage levels, and time points, along with an action set of water outflows. State transitions are deduced from the inflow distributions. A cost is determined for each state-action pair from the cost of fuel and curtailment required to meet system load net of hydropower generation.

An optimal operating policy and corresponding offer curves were computed from a pair of dual LPs. We presented quantitative results using inflow data from a New Zealand reservoir spanning seven decades.

Potential extensions to the described approach include:
\begin{enumerate}
    \item Separately modeling the flows of multiple reservoirs.
    \item Replacing the (risk-neutral) expected cost in LP \eqref{eqn:state-action-lp} with a risk-averse objective (e.g., conditional value-at-risk \cite{Rockafellar-Uryasev:2000:CVaR}).
    \item Including additional performance constraints (e.g., to put probabilistic limits on the incidence of demand curtailment).
    \item Using higher-fidelity problem data (e.g., time-varying system load and fuel prices).
\end{enumerate}

\section*{Acknowledgment}

The authors thank Prof. Benjamin Hobbs (Johns Hopkins University, Whiting School of Engineering) for valuable discussions.

This material is based upon work supported by the U.S. Department of Energy's Office of Energy Efficiency and Renewable Energy (EERE) under the Wind Energy Technologies Office (WETO) Award Number DE-EE0011269, the Massachusetts Clean Energy Center and the Maryland Energy Administration. The views expressed herein do not necessarily represent the views of the U.S. Department of Energy, the United States Government, the Massachusetts Clean Energy Center or the Maryland Energy Administration.

\bibliographystyle{IEEEtran}
\bibliography{references.bib}

\end{document}